 \documentclass[11pt]{article}
 \usepackage{amsfonts,amsmath,amssymb,mathrsfs}
 \usepackage{color}
 \usepackage[svgnames]{xcolor}
 \usepackage{harvard}
 \usepackage{bm}

 \allowdisplaybreaks

 \newtheorem{theorem}{Theorem}[section]
 \newtheorem{lemma}[theorem]{Lemma}
 \newtheorem{corollary}[theorem]{Corollary}
 \newtheorem{proposition}[theorem]{Proposition}
 \newtheorem{remark}[theorem]{Remark}
 \newtheorem{example}[theorem]{Example}
 \newtheorem{condition}[theorem]{Condition}
 \newtheorem{definition}[theorem]{Definition}

 \def\blemma{\begin{lemma}\sl{}
 \def\elemma{\end{lemma}}}
 
 \def\bproposition{\begin{proposition}\sl{}
 \def\eproposition{\end{proposition}}}

 \def\benumerate{\begin{enumerate}}\def\eenumerate{\end{enumerate}}
 \def\bitemize{\begin{itemize}}\def\eitemize{\end{itemize}}

 \def\beqlb{\begin{eqnarray}}
 \def\eeqlb{\end{eqnarray}}
 \def\beqnn{\begin{eqnarray*}}
 \def\eeqnn{\end{eqnarray*}}

 \def\qed{\hfill$\Box$\medskip}
 \def\bproof{\begin{proof}}\def\eproof{\qed\end{proof}}

 \def\<{\langle}\def\>{\rangle}

 \def\mbb{\mathbb}
 \def\mbf{\mathbf}\def\mrm{\mathrm}

 \def\ar{\!\!&}

 \def\d{\mrm{d}}\def\e{\mrm{e}}
 \def\sgn{\mathrm{sgn}}

\begin{document}

$~$

\bigskip

\centerline{\Large\bf Two-type continuous-state branching processes}

\smallskip

\centerline{\Large\bf in varying environments\footnote{The research is supported by the National Key R\&D Program of China (No. 2020YFA0712901).}}

\bigskip

\centerline{Zenghu Li\footnote{E-mail: lizh@bnu.edu.cn. Laboratory of Mathematics and Complex Systems, School of Mathematical Sciences, Beijing Normal University, Beijing 100875, People's Republic of China.}
 \quad
Junyan Zhang$^{*}$\footnote{*corresponding author, E-mail: zhangjunyan@mail.bnu.edu.cn. Laboratory of Mathematics and Complex Systems, School of Mathematical Sciences, Beijing Normal University, Beijing 100875, People's Republic of China.}}

 \begin{center}
 \begin{minipage}{12cm}
 \begin{center}\textbf{Abstract}\end{center}
 \footnotesize

A basic class of two-type continuous-state branching processes in varying environments are constructed by solving the backward equation determining the cumulant semigroup. The parameters of the process are allowed to be c\`adl\`ag in time and the difficulty brought about by the bottlenecks are overcome by introducing a suitable moment condition.

 \end{minipage}
 \end{center}

\bigskip\bigskip



\section{Introduction}

\setcounter{equation}{0}

Continuous-state branching processes in varying environments (CBVE-processes) are positive inhomogeneous Markov processes arising as scaling limits of discrete Galton--Watson branching processes. A limit theorem of this type for one-dimensional processes was established by Bansaye and Simatos (2015), which generalizes the results for the homogeneous processes; see, e.g., Aliev and Shchurenkov (1982), Grimvall (1974), Lamperti (1967) and Li (2020). The distributional properties of a CBVE-process is determined by its \textit{cumulant semigroup}, which gives the exponents of the Laplace transforms of its transition probabilities. The semigroup was defined by a backward integral equation involving the \textit{branching mechanism} of the process. In the inhomogeneous case, the uniqueness of the solution to the backward equation is an annoying problem because of the possible presence of the \textit{bottlenecks}, which are times when the solution jumps to zero. In fact, the determination of the behavior of a general CBVE-process at the bottlenecks was left open in Bansaye and Simatos (2015). The problem was settled in Fang and Li (2022), where the CBVE-process was constructed as the pathwise unique solution to a stochastic integral equation. The study of inhomogeneous cumulant semigroups is closely related with that of reverse evolution families in complex analysis; see Gumenyuk et al. (2024, 2022+).

For homogeneous multi-dimensional continuous-state branching processes, a uniqueness problem for the backward equation of the cumulant semigroup was pointed out in Rhyzhov and Skorokhod (1970), who expected that the uniqueness of solution holds for all initial values if and only if it holds for the initial value zero. A proof of their assertion was given in the very recent work of Li and Li (2024).

In this note, we discuss the construction of \textit{two-dimensional continuous-state branching processes in varying environments} (TCBVE-processes) under a suitable moment condition. By a TCBVE-process we mean an inhomogeneous Markov process $\mbf{X}= \{(X_1(t),X_2(t)): t\ge 0\}$ in the state space $\mbb{R}_+^2$ with transition semigroup $(Q_{r,t})_{t\geq r\ge 0}$ defined by
 \beqlb\label{Laplace transform}
\int_{\mbb{R}_+^2}\e^{-\langle {\bm{\lambda}},\mbf{y}\rangle}Q_{r,t}(\mbf{x},\d\mbf{y})
 =
\e^{-\langle \mbf{x},\mbf{v}_{r,t}({\bm\lambda})\rangle},\quad\bm\lambda,\mbf{x}\in \mbb{R}_+^2,
 \eeqlb
where $\langle\cdot, \cdot\rangle$ denotes the Euclidean inner product and $(\mbf{v}_{r,t})_{t\geq r\ge 0}= (v_{1,r,t},v_{2,r,t})_{t\geq r\ge 0}$ is a family of continuous transformations on $\mbb{R}_+^2$. From \eqref{Laplace transform} it is easy to see that each $Q_{r,t}(\mbf{x},\cdot)$ is an infinitely divisible probability measure on $\mbb{R}_+^2$. Then we have the \textit{L\'{e}vy--Khintchine type representation}, for $i=1,2$,
 \beqlb\label{LK representation}
v_{i,r,t}(\bm\lambda)= \langle\mbf{h}_{i,r,t},\bm\lambda\rangle + \int_{\mbb{R}_+^2\setminus\{\mbf{0}\}}(1-\e^{-\langle\bm\lambda,\mbf{y}\rangle})l_{i,r,t}(\d\mbf{y}),
 \eeqlb
where $\mbf{h}_{i,r,t}\in\mbb{R}_+^2$ and $(1\wedge\Vert\mbf{y}\Vert)l_{i,r,t}(\d\mbf{y})$ is a finite measure on $\mbb{R}_+^2\setminus\{\mbf{0}\}$. The Chapman--Kolmogorov equation of $(Q_{r,t})_{t\geq r\ge 0}$ implies
 \beqlb\label{semigroup property}
\mbf{v}_{r,t}(\bm\lambda)=\mbf{v}_{r,s}\circ \mbf{v}_{s,t}(\bm\lambda), \qquad \bm\lambda\in \mbb{R}_+^2,~ t\geq s\geq r\ge 0.
 \eeqlb
We call $(\mbf{v}_{r,t}(\cdot))_{t\geq r\ge 0}$ the \textit{cumulant semigroup} of the TCBVE-process $\mbf{X}$.
	
For $i=1,2$, let $c_i$ be an increasing continuous function on $[0,\infty)$ satisfying $c_i(0)=0$ and let $b_{ii}$ be a c\`adl\`ag function on $[0,\infty)$ satisfying $b_{ii}(0)=0$ and having locally bounded variations. Let $m_i$ be a $\sigma$-finite measure on $(0,\infty)\times (\mbb{R}_+^2\setminus\{\mbf{0}\})$ satisfying
 \beqlb\label{m_i(t)}
m_i(t):= \int_0^t\int_{\mbb{R}_+^2\setminus\{\mbf{0}\}}(z_i^2\mbf 1_{\{\Vert\mbf{z}\Vert\leq 1\}}+z_i\mbf 1_{\{\Vert\mbf{z}\Vert> 1\}}+z_j)m_i(\d s,\d\mbf{z})<\infty, \quad t\geq 0.
 \eeqlb
Moreover, we assume that
 \beqlb\label{delta_i}
\delta_i(t)= \Delta b_{ii}(t) + \int_{\mbb{R}_+^2\setminus\{\mbf{0}\}}z_im_i(\{t\},\d\mbf{z})\leq 1, \quad t\geq 0,
 \eeqlb
where $\Delta b_{ii}(t)= b_{ii}(t) - b_{ii}(t-)$. For $i,j=1,2$ with $i\neq j$, let $b_{ij}$ be an increasing c\`adl\`ag function on $[0,\infty)$ satisfying $b_{ij}(0)=0$. For $\bm\lambda\in\mbb{R}_+^2$ and $t\geq 0$, we consider the following system of backward integral equations
 \beqlb\label{backward eq}
v_{i,r,t}(\bm\lambda)\ar=\ar \lambda_i + \int_r^t v_{j,s,t}(\bm\lambda)b_{ij}(\d s) - \int_r^t v_{i,s,t}(\bm\lambda)b_{ii}(\d s) - \int_r^t v_{i,s,t}(\bm\lambda)^2c_i(\d s) \cr
 \ar\ar
- \int_r^t\int_{\mbb{R}_+^2\setminus\{\mbf{0}\}}K_i(\mbf{v}_{s,t}(\bm\lambda),\mbf{z})m_i(\d s,\d\mbf{z}), \qquad r\in [0,t],~ i=1,2,
 \eeqlb
where
 $$
K_i(\bm\lambda,\mbf{z})=\e^{-\langle \bm\lambda,\mbf{z}\rangle}-1+\lambda_i z_i.
 $$
Throughout this note, we make the conventions
 \beqnn
\int_r^t=\int_{(r,t]}, ~ \int_r^\infty=\int_{(r,\infty)}, \qquad t\geq r\in \mbb{R}.
 \eeqnn
The main results of the note are the following:

\begin{theorem}\label{main result} Under the assumptions described above, for any $t\geq 0$ and $\bm\lambda\in \mbb{R}_+^2$ there is a unique bounded solution $[0,t]\ni r\mapsto \mbf{v}_{r,t}(\bm\lambda)\in \mbb{R}_+^2$ to \eqref{backward eq}. Moreover, an inhomogeneous transition semigroup $(Q_{r,t})_{t\geq r\ge 0}$ on $\mbb{R}_+^2$ is defined by \eqref{Laplace transform}. \end{theorem}

\begin{theorem}\label{main result(2)} Let $(Q_{r,t})_{t\geq r\ge 0}$ be the inhomogeneous transition semigroup on $\mbb{R}_+^2$ defined by \eqref{Laplace transform}. Let
 \beqlb\label{b_ij-}
\bar{b}_{ij}(t)=b_{ij}(t)+\int_0^t\int_{\mbb{R}_+^2\setminus\{\mbf{0}\}}z_jm_i(\d s,\d\mbf{z}), \quad t\ge 0.
 \eeqlb
Then we have
 \beqlb\label{first moments}
\int_{\mbb{R}_+^2} \langle {\bm{\lambda}},\mbf{y}\rangle Q_{r,t}(\mbf{x},\d\mbf{y})
 =
\langle \mbf{x},\bm{\pi}_{r,t}({\bm\lambda})\rangle,\qquad \bm\lambda\in \mbb{R}^2,~ \mbf{x}\in \mbb{R}_+^2,
 \eeqlb
where $[0,t]\ni r\mapsto \bm{\pi}_{r,t}(\bm\lambda)\in \mbb{R}^2$ is the unique bounded solution to
 \beqlb\label{linear backward eq}
\pi_{i,r,t}(\bm\lambda)= \lambda_i + \int_r^t \pi_{j,s,t}(\bm\lambda)\bar{b}_{ij}(\d s) - \int_r^t \pi_{i,s,t}(\bm\lambda)b_{ii}(\d s), \quad r\in [0,t].
 \eeqlb
\end{theorem}

The assumption \eqref{delta_i} is necessary to guarantee that $\mbf{v}_{r,t}(\bm\lambda)\in \mbb{R}_+^2$ for every $r\in [0,t]$ and $\bm\lambda\in \mbb{R}_+^2$. The detailed explanation of a similar assumption for the one-dimensional CBVE-process was given in Fang and Li (2022, Remark~1.7). The moment condition \eqref{m_i(t)} is slightly stronger than those in Bansaye and Simatos (2015) and Fang and Li (2022). Under this condition, the functions
 \beqlb\label{lam-to-intKdm} 
\bm\lambda\mapsto \int_r^t\int_{\mbb{R}_+^2\setminus\{\mbf{0}\}} K_i(\bm\lambda,\mbf{z})m_i(\d s,\d\mbf{z}), \quad i=1,2,\ t\geq r\geq 0
 \eeqlb
are Lipschitz in each bounded subset of $\mbb{R}_+^2$. This \textit{local-Lipschitz property} is necessary to guarantee the uniqueness of the solution to \eqref{backward eq}. Without condition \eqref{m_i(t)}, the above functions could behave rather irregularly near the boundary $\partial \mbb{R}_+^2:= (\mbb{R}_+\times \{0\})\cup (\{0\}\times \mbb{R}_+)$. For $i=1,2$ let
 \beqnn
K_i= \{s> 0: \Delta b_{ii}(s)= 1, \Delta b_{ij}(s)= m_i(\{s\}\times (\mbb{R}_+^2\setminus\{\mbf{0}\}))= 0\}.
 \eeqnn
We call $K:= K_1\cup K_2$ the set of \textit{bottlenecks}. It is easy to see that $K\cap (0,t]$ is a finite set for each $t\ge 0$. Let $p(t)= \max K\cap (0,t]$ be the \textit{last bottleneck} before time $t\ge 0$. If $K\cap (0,t]\neq \emptyset$, then $0< p(t)\le t$ and one can see from \eqref{backward eq} that either $v_{1,p(t)-,t}(\bm\lambda)= 0$ or $v_{2,p(t)-,t}(\bm\lambda)= 0$, and hence $\mbf{v}_{p(t)-,t}(\bm\lambda)\in \partial \mbb{R}_+^2$. Therefore, the bottlenecks and the irregularities of the functions in \eqref{lam-to-intKdm} near $\partial \mbb{R}_+^2$ would bring difficulties to the discussion of the equation \eqref{backward eq} if condition \eqref{m_i(t)} were not assumed. The above theorems give characterizations of a basic class of TCBVE-processes and serves as the basis of constructions of more general processes, which are carried out in a forthcoming work.

The rest of the note is organized as follows. In Section~2, we study a special form of the equation \eqref{backward eq}. The proofs of Theorems~\ref{main result} and \ref{main result(2)} are given in Section~3.

\section{Preliminaries}

\setcounter{equation}{0}

In this section, we study a special kind of backward equation systems, which can be obtained by iteration and inhomogeneous nonlinear h-transformation. We first consider the two-dimensional Gronwall's inequality, which is essential in the proof of uniqueness of the solution to backward equation system.

\blemma\label{multiGronwall}
Let $\beta_{ij}:[0,T]\to\mbb{R}$, $i,j=1,2$ be a right-continuous nondecreasing function. Let $a_i:[0,T]\to\mbb{R}_+$, $i=1,2$ be a right-continuous nondecreasing function. Let $g_i:[0,T]\to\mbb{R}_+$, $i=1,2$ be a measurable function such that
\beqnn
	\int_0^T[g_1(s)+g_2(s)]\d\beta(s)<\infty,
\eeqnn
where $\beta(t)=\beta_{11}(t)+\beta_{12}(t)+\beta_{21}(t)+\beta_{22}(t)$. Suppose for all $t\in[0,T]$ and $(i,j)=(1,2)\ \text{or}\ (2,1)$, we have
\beqnn
	g_i(t)\leq a_i(t)+\int_0^t g_i(s)\d\beta_{ii}(s)+\int_0^t g_j(s)\d\beta_{ij}(s).
\eeqnn
Then
\beqlb\label{multiGronwallineq}
	g_i(t)\leq d_i(t)\exp\left\{\int_0^t\int_s^t \e^{\beta_{jj}(r)-\beta_{jj}(0)}\d\beta_{ij}(r)\d\beta_{ji}(s)+\beta_{ii}(t)-\beta_{ii}(0)\right\},
\eeqlb
where
\beqnn
	d_i(t)=a_i(t)+\int_0^t a_j(s)\e^{\beta_{jj}(s)-\beta_{jj}(0)}\d\beta_{ij}(s).
\eeqnn
\elemma

\bproof
By Gronwall's inequality, see Lemma 2.1 in Mao (1990), we have
\beqnn
	g_j(t)\leq \left[a_j(t)+\int_0^t g_i(s)\d\beta_{ji}(s)\right]\exp\{\beta_{jj}(t)-\beta_{jj}(0)\},
\eeqnn
then
\beqnn
	g_i(t)\leq a_i(t)+\int_0^t a_j(s)\e^{\beta_{jj}(s)-\beta_{jj}(0)}\d\beta_{ij}(s)+\int_0^t g_i(s)\d\beta_{ii}(s)\\
	+\int_0^t g_i(s)\int_s^t \e^{\beta_{jj}(r)-\beta_{jj}(0)}\d\beta_{ij}(r)\d\beta_{ji}(s).
\eeqnn
We can easily get \eqref{multiGronwallineq} by Gronwall's inequality.
\eproof

Let $\gamma_{ii}, \gamma_{ij}$ be c\`adl\`ag functions on $[0,\infty)$ with locally bounded variations. Suppose that $\Delta\gamma_{ii}(t)>-1$ holds for all $t>0$ and $t\mapsto \gamma_{ij}(t)$ is increasing. Let $(z_1+z_2)\mu_i(\d s,\d\mbf{z})$ be a finite measure on $(0,\infty)\times \big(\mbb{R}_+^2\setminus\{\mbf{0}\}\big)$. The following Lemma is a special two-dimensional inhomogeneous nonlinear h-transformation.

\begin{lemma}\label{2.3 h-transform}
	Suppose that $((\bm u_{r,t}(\bm\lambda))_{t\geq r}$ is a positive solution to
	\beqlb\label{2.3 backward eq}
	\begin{split}
		u_{i,r,t}(\bm\lambda)=\lambda_i &+\int_r^t u_{i,s,t}(\bm\lambda)\gamma_{ii}(\d s)+\int_r^t u_{j,s,t}(\bm\lambda)\gamma_{ij}(\d s)\\
		&+\int_r^t\int_{\mbb{R}_+^2\setminus\{\mbf{0}\}}\left(1-\e^{-\langle\mbf u_{s,t}(\bm\lambda),\mbf{z}\rangle}\right)\mu_i(\d s,\d\mbf{z}),\quad r\in[0,t].
	\end{split}
\eeqlb
	Let $t\mapsto \zeta_i(t)$ be c\`adl\`ag function on $[0,\infty)$ with locally bounded variations. Define
	\beqlb\label{h-transform semigroup}
		 v_{i,r,t}(\bm\lambda)=\e^{\zeta_i(r)}u_{i,r,t}(\e^{-\zeta_1(t)}\lambda_1,\e^{-\zeta_2(t)}\lambda_2),\quad\bm\lambda\in\mbb{R}_+^2.
	\eeqlb
	Then $r\mapsto \bm v_{r,t}(\bm\lambda)$ is a positive solution to
	\beqlb\label{2.3 h-transform eq}
	\begin{split}
		v_{i,r,t}(\bm\lambda)&=\lambda_i-\int_r^t v_{i,s,t}(\bm\lambda)\eta_i(\d s)+\int_r^t \e^{-\Delta\zeta_i(s)}v_{i,s,t}(\bm\lambda)\gamma_{ii}(\d s)\\
		&\quad+\int_r^t \e^{\zeta_i(s-)-\zeta_j(s)}v_{j,s,t}(\bm\lambda)\gamma_{ij}(\d s)\\
		&\quad+\int_r^t\int_{\mbb{R}_+^2\setminus\{\bm 0\}}\left(1-\e^{-\langle\bm v_{s,t}(\bm\lambda),\bm z\rangle}\right)\e^{\zeta_i(s-)}\mu_i(\d s,\e^{\zeta_1(s)}\d z_1,\e^{\zeta_2(s)}\d z_2),
	\end{split}	
	\eeqlb
	where
	\beqnn
		\eta_i(s)=\zeta_{i,c}(s)+\sum_{s\in(0,t]}\left(1-\e^{-\Delta\zeta_i(s)}\right).
	\eeqnn

\end{lemma}
\bproof
	By integration by parts we have
	\beqnn
		\lambda_i\ar=\ar \e^{\zeta_i(r)}u_{i,r,t}(\e^{-\zeta_1(t)}\lambda_1,\e^{-\zeta_2(t)}\lambda_2)+\int_r^t u_{i,s,t}(\e^{-\zeta_1(t)}\lambda_1,\e^{-\zeta_2(t)}\lambda_2)\d \e^{\zeta_i(s)}\\
		\ar\ar+\int_r^t \e^{\zeta_i(s-)}\d u_{i,s,t}(\e^{-\zeta_1(t)}
		\lambda_1,\e^{-\zeta_2(t)}\lambda_2)\\
		\ar=\ar v_{i,r,t}(\bm\lambda)+\int_r^t u_{i,s,t}(\e^{-\zeta_1(t)}\lambda_1,\e^{-\zeta_2(t)}\lambda_2)\e^{\zeta_i(s)}\zeta_{i,c}(\d s)\\
		\ar\ar+\sum_{s\in(r,t]}u_{i,s,t}(\e^{-\zeta_1(t)}\lambda_1,\e^{-\zeta_2(t)}\lambda_2)\left(\e^{\zeta_i(s)}-\e^{\zeta_i(s-)}\right)\\
		\ar\ar-\int_r^t \e^{\zeta_i(s-)} u_{i,s,t}(\e^{-\zeta_1(t)}\lambda_1,\e^{-\zeta_2(t)}\lambda_2)\gamma_{ii}(\d s)\\
		\ar\ar-\int_r^t \e^{\zeta_i(s-)}u_{j,s,t}(\e^{-\zeta_1(t)}\lambda_1,\e^{-\zeta_2(t)}\lambda_2)\gamma_{ij}(\d s)\\
		\ar\ar-\int_r^t\int_{\mbb{R}_+^2\setminus\{\bm 0\}}\left(1-\e^{-\langle\bm u_{s,t}(\e^{-\zeta_1(t)}\lambda_1,\e^{-\zeta_2(t)}\lambda_2),\bm z\rangle}\right)\e^{\zeta_i(s-)}\mu_i(\d s,\d\mbf{z})\\
		\ar=\ar v_{i,r,t}(\bm\lambda)+\int_r^t v_{i,s,t}(\bm\lambda)\zeta_{i,c}(\d s)+\sum_{s\in(r,t]}v_{i,s,t}(\bm\lambda)\left(1-\e^{-\Delta\zeta_i(s)}\right)\\
		\ar\ar-\int_r^t \e^{-\Delta\zeta_i(s)} v_{i,s,t}(\bm\lambda)\gamma_{ii}(\d s)-\int_r^t \e^{\zeta_i(s-)-\zeta_j(s)}v_{j,s,t}(\bm\lambda)\gamma_{ij}(\d s)\\
		\ar\ar-\int_r^t\int_{\mbb{R}_+^2\setminus\{\bm 0\}}\left(1-\e^{-\langle\bm v_{s,t}(\bm\lambda),\bm z\rangle}\right)\e^{\zeta_i(s-)}\mu_i(\d s,\e^{\zeta_1(s)}\d z_1,\e^{\zeta_2(s)}\d z_2).
	\eeqnn
\eproof

\begin{remark}
	\eqref{2.3 backward eq} is a special case of \eqref{backward eq} with $c_i\equiv 0$, $b_{ij}(\d s)=\gamma_{ij}(\d s)$, $m_i(\d s,\d\mbf{z})=\mu_i(\d s,\d\mbf{z})$ and
	\beqnn
		b_{ii}(\d s)=-\gamma_{ii}(\d s)-\int_{\mbb{R}_+^2\setminus\{\mbf{0}\}}z_i\mu_i(\d s,\d\mbf{z}).
	\eeqnn
\end{remark}

Now let us consider the uniqueness and existence of the solution to backward equation system \eqref{2.3 backward eq}.
\begin{proposition}\label{special case backward eq solution}
For $t\geq 0$ and $\bm\lambda\in\mbb{R}_+^2$, there is a unique bounded positive solution $r\mapsto \mbf u_{r,t}(\bm\lambda)$ on $[0,t]$ to \eqref{2.3 backward eq} and $(\mbf u_{r,t}(\bm\lambda))_{t\geq r}$ is a cumulant semigroup. Moreover, for $\bm\lambda\in\mbb{R}_+^2$, we have
\beqlb\label{2.3 estimate}
	u_{1,r,t}(\bm\lambda)+u_{2,r,t}(\bm\lambda)\leq (\lambda_1+\lambda_2)\e^{\rho(r,t]},
\eeqlb
where
\beqnn
	\rho(t)=\Vert\rho_1\Vert(t)+\Vert\rho_2\Vert(t)+\gamma_{12}(t)+\gamma_{21}(t)+\int_0^t\int_{\mbb{R}_+^2\setminus\{\mbf{0}\}}z_1\mu_2(\d s,\d\mbf{z})+\int_0^t\int_{\mbb{R}_+^2\setminus\{\mbf{0}\}}z_2\mu_1(\d s,\d\mbf{z})
\eeqnn	
and
\beqlb\label{rho_i}
	 \rho_i(t)=\gamma_{ii}(t)+\int_0^t\int_{\mbb{R}_+^2\setminus\{\mbf{0}\}}z_i\mu_i(\d s,\d\mbf{z}).
\eeqlb
\end{proposition}

\bproof
	Step 1. Let $r\mapsto \mbf u_{r,t}(\bm\lambda)$ be a positive solution to \eqref{2.3 backward eq}. By elementary calculus, we can easily obtain
	\beqnn
	\begin{aligned}
		u_{1,r,t}(\bm\lambda)&+u_{2,r,t}(\bm\lambda)\leq\int_r^t u_{1,s,t}(\bm\lambda)\left(\Vert \rho_1\Vert(\d s)+\gamma_{21}(\d s)+\int_{\mbb{R}_+^2\setminus\{\mbf{0}\}}z_1\mu_2(\d s,\d\mbf{z})\right)\\
		&+\lambda_1+\lambda_2+\int_r^t u_{2,s,t}(\bm\lambda)\left(\Vert \rho_2\Vert(\d s)+\gamma_{12}(\d s)+\int_{\mbb{R}_+^2\setminus\{\mbf{0}\}}z_2\mu_1(\d s,\d\mbf{z})\right),
	\end{aligned}		
	\eeqnn
	where $\rho_i(t)$ is defined by \eqref{rho_i}. Then we have the upper bound estimation \eqref{2.3 estimate} by Gronwall's inequality. Suppose that $r\mapsto \mbf w(r,t,\bm\lambda)$ is also a positive solution to \eqref{2.3 backward eq}. Then we have
	\beqnn
		|u_{i,r,t}(\bm\lambda)-w_{i,r,t}(\bm\lambda)|\leq \int_r^t[|u_{i,s,t}(\bm\lambda)-w_{i,s,t}(\bm\lambda)|+|u_{j,s,t}(\bm\lambda)-w_{j,s,t}(\bm\lambda)|]\Vert\rho\Vert(\d s).
	\eeqnn
	By Lemma \ref{multiGronwall}, we see $|u_{i,r,t}(\bm\lambda)-w_{i,r,t}(\bm\lambda)|=0$ for all $r\in[0,t]$, implying the uniqueness. Note that the uniqueness of the solution implies the semigroup property.
	
	Step 2. We only need to prove the case when $\gamma_{ii}$ vanishes. In fact, suppose that $r\mapsto \mbf{v}_{r,t}(\bm\lambda)$ is a solution to
	\beqlb
	\begin{split}
		v_{i,r,t}(\bm\lambda)=\lambda_i &+\int_r^t\int_{\mbb{R}_+^2\setminus\{\mbf{0}\}}\big(1-\e^{-\langle \mbf{v}_{s,t}(\bm\lambda),\mbf{z}\rangle}\big)\e^{\zeta_i(s-)}\mu_i(\d s,\e^{\zeta_1(s)}\d z_1,\e^{\zeta_2(s)}\d z_2)\\
		&+\int_r^t \e^{\zeta_i(s-)-\zeta_j(s)}v_{j,s,t}(\bm\lambda)\gamma_{ij}(\d s),
	\end{split}
	\eeqlb
	where $\zeta_i$ is the c\`adl\`ag function on $[0,\infty)$ such that $\zeta_{i,c}(t)=\gamma_{ii,c}(t)$ and $\Delta\zeta_i(t)=\log[1+\Delta\gamma_{ii}(t)]$ for every $t>0$. Define $u_{i,r,t}(\bm\lambda)=\e^{-\zeta_i(r)}v_{i,r,t}(\e^{\zeta_1(t)}\lambda_1, \e^{\zeta_2(t)}\lambda_2)$. Then by Lemma \ref{2.3 h-transform}, we have $r\mapsto\mbf u_{r,t}(\bm\lambda)$ is a solution to \eqref{2.3 backward eq}. Moreover, if $(\mbf{v}_{r,t}(\bm\lambda))_{t\geq r}$ is a cumulant semigroup, then $(\mbf u_{r,t}(\bm\lambda))_{t\geq r}$ is also a cumulant semigroup by the uniqueness and L\'evy-Kthintchine representation of the solution.
		
	Step 3. Given $t\geq 0$ and $\bm\lambda\in\mbb{R}_+^2$, let $v_{i,r,t}^{(0)}(\bm\lambda)=\lambda_i$ and define $v_{i,r,t}^{(k)}(\bm\lambda)$ inductively by
	\beqnn
		v_{i,r,t}^{(k+1)}(\bm\lambda)=\lambda_i+\int_r^t\int_{\mbb{R}_+^2\setminus\{\mbf{0}\}}\big(1-\e^{-\langle \mbf{v}^{(k)}_{s,t}(\bm\lambda),\mbf{z}\rangle}\big)\mu_i(\d s,\d\mbf{z})+\int_r^t v_{j,s,t}^{(k)}(\bm\lambda)\gamma_{ij}(\d s),
	\eeqnn
	where $r\in[0,t]$. By Watanabe (1969) one can easily see that $v_{i,r,t}^{(k)}(\bm\lambda)$ has the L\'evy-Kthintchine representation \eqref{LK representation}, inductively. Furthermore, we can easily obtain
 	\beqnn
		v_{i,r,t}^{(k)}(\mbf{0})=0\leq v_{i,r,t}^{(k)}(\bm\lambda)\leq v_{i,r,t}^{(k+1)}(\bm\lambda)\leq 2\Vert\bm\lambda\Vert \e^{\rho(t)}.
	\eeqnn
	Let
	\beqnn
		u(k,r,t,\bm\lambda)=\sup_i\sup_{r\leq s\leq t}|v_{i,s,t}^{(k)}(\bm\lambda)-v_{i,s,t}^{(k-1)}(\bm\lambda)|,\quad k\geq 1
	\eeqnn
	and $u(0,r,t,\bm\lambda):=\lambda_1\vee\lambda_2$. Then, for every $\bm\lambda\in[0,B]^2$, we have
	\beqnn
		u(k,r,t,\bm\lambda)\ar\leq\ar 2\int_r^t\int_{\mbb{R}_+^2\setminus\{\mbf{0}\}}u(k-1,t_1,t,\bm\lambda)\rho(\d t_1)\\
			\ar\leq\ar 2^2\int_r^t \rho(\d t_1)\int_{t_1}^t u(k-2,t_2,t,\bm\lambda)\rho(\d t_2)\leq \cdots\\
			\ar\leq\ar 2^{k}\int_r^t \rho(\d t_1)\int_{t_1}^t\cdots\int_{t_{k-1}}^t u(0,t_{k},t,\bm\lambda)\rho(\d t_{k})\\
			\ar\leq\ar 2^{k}B\int_r^t \rho(\d t_1)\int_{t_1}^t\cdots\int_{t_{k-1}}^t \rho(\d t_{k})\\
			\ar\leq\ar B\frac{(2\rho(t))^{k}}{k!},
	\eeqnn
	and hence $\sum_{k=1}^\infty u(k,r,t,\bm\lambda)\leq B\e^{2\rho(t)}<\infty$. Thus the limit $v_{i,r,t}(\bm\lambda):=\uparrow\lim_{k\to\infty}v_{i,r,t}^{(k)}(\bm\lambda)$ exists and the convergence is uniform in $(r,\bm\lambda)\in [0,t]\times[0,B]\times[0,B]$ for every $t\geq 0$ and $B\geq 0$. By monotone convergence theorem, $r\mapsto \mbf{v}_{r,t}(\bm\lambda)$ is a solution to \eqref{2.3 backward eq} with $\gamma_{ii}\equiv 0$. Moreover, the limit $v_{i,r,t}(\bm\lambda)$ also has the L\'evy-Kthintchine representation \eqref{LK representation}, see Lemma 1 in Watanabe (1969), implying $(\mbf{v}_{r,t}(\bm\lambda))_{t\geq r}$ is a cumulant semigroup.
\eproof

\section{General backward equations}

\setcounter{equation}{0}

In this section, we use solutions to special backward equation systems \eqref{2.3 backward eq} to approximate the solution to \eqref{backward eq}, which is similar to Fang and Li (2022). For simplicity, we use another way to describe the backward equation system \eqref{backward eq}. Let $B[0,\infty)^+$ be the set of locally bounded positive Borel functions on $[0,\infty)$ and $\mbf B[0,\infty)^+$ $:=\{\mbf f=(f_1,f_2):f_1,f_2\in B[0,\infty)^+\}$. Let $\mathscr B[0,\infty)$ be the set of Borel sets on $[0,\infty)$.

\begin{definition} Let $\bar b_{ij},b_{ii},c_i,m_i$ be defined as in introduction. For $\mbf f\in \mbf B[0,\infty)^+$, $B\in\mathscr B[0,\infty)$, let $\phi_i$ be a functional on $\mbf B[0,\infty)^+ \times \mathscr B[0,\infty)$ defined by
\beqlb\label{phi}
\begin{split}
	\phi_i(\mbf f,B)=\int_B f_i(s)b_{ii}(\d s)&-\int_B f_j(s) \bar{b}_{ij}(\d s)+\int_B f_i^2(s)c_i(\d s)\\
	&+\int_B\int_{\mbb{R}_+^2\setminus\{\mbf{0}\}}K(\mbf f(s),\mbf{z})m_i(\d s,\d\mbf{z}),
\end{split}
\eeqlb
where $K(\bm\lambda,\mbf{z})=\e^{-\langle\bm\lambda,\mbf{z}\rangle}-1+\langle\bm\lambda,\mbf{z}\rangle$. We call $\bm\phi=(\phi_1,\phi_2)$ a \textit{branching mechanism} with parameters $(\bar{b}_{ij},b_{ii},c_i,m_i)$. \end{definition}

We denote $\mbf e_i$ the vector that $i$-th component is $1$ and the remaining component is $0$. For $n\geq 1$, define a branching mechanism $\bm\phi_n$ by
\beqlb\label{phi_n}\small
	\begin{split}
		\phi_{n,i}(\mbf f,B)&=\int_B f_i(s)b_{ii}(\d s)-\e^{-n}\int_B f_i(s)\Vert b_{ii}\Vert(\d s)-\int_B f_j(s)\bar{b}_{ij}(\d s) \\
		&\quad+2n^2\int_B \big(\e^{-f_i(s)/n}-1+f_i(s)/n\big)c_i(\d s)\\
		&\quad+(1-\e^{-n})\int_B\int_{\mbb{R}_+^2\setminus\{\mbf{0}\}}K (\mbf f(s),\mbf{z})(1\wedge n\Vert \mbf{z}\Vert)m_i(\d s,\d\mbf{z})\\
		&=-\int_B f_i(s)\gamma_{n,ii}(\d s)-\int_B f_j(s)\gamma_{n,ij}(\d s)-\int_B\int_{\mbb{R}_+^2\setminus\{\mbf{0}\}}(1-\e^{-\langle \mbf f(s),\mbf{z}\rangle})\mu_{n,i}(\d s,\d\mbf{z}),
	\end{split}
\eeqlb
where
\beqnn
	\gamma_{n,ii}(\d s)\ar=\ar-b_{ii}(\d s)+\e^{-n}\Vert b_{ii}\Vert(\d s)-2nc_i(\d s)-(1-\e^{-n})\int_{\mbb{R}_+^2\setminus\{\mbf{0}\}}z_i(1\wedge n\Vert \mbf{z}\Vert)m_i(\d s,\d\mbf{z}),\\
	\gamma_{n,ij}(\d s)\ar=\ar\bar{b}_{ij}(\d s)-(1-\e^{-n})\int_{\mbb{R}_+^2\setminus\{\mbf{0}\}}z_j(1\wedge n\Vert \mbf{z}\Vert)m_i(\d s,\d\mbf{z})
\eeqnn
and
\beqnn\small
		\mu_{n,i}(\d s,\d\mbf{z})\!=\!2n^2c_i(\d s)\delta_{\mbf e_i}(n\d\mbf{z})+(1-\e^{-n})(1\wedge n\Vert \mbf{z}\Vert)m_i(\d s,\d\mbf{z}).
\eeqnn
Then it is easy to see that $\Delta\gamma_{n,ii}(t)>-1$ and $\gamma_{n,ij}(\d s)$ is a measure on $[0,\infty)$. In fact, we have
\beqnn
	\Delta\gamma_{n,ii}(t)\ar=\ar-\Delta b_{ii}(t)+\e^{-n}\Delta\Vert b_{ii}\Vert(t)-(1-\e^{-n})\int_{\mbb{R}_+^2\setminus\{\mbf{0}\}}z_i(1\wedge n\Vert \mbf{z}\Vert)m_i(\{t\},\d\mbf{z})\\
	\ar\geq\ar-\delta_i(t)+\e^{-n}\Delta\Vert b_{ii}\Vert(t)+\e^{-n}\int_{\mbb{R}_+^2\setminus\{\mbf{0}\}}z_i(1\wedge n\Vert \mbf{z}\Vert)m_i(\{t\},\d\mbf{z})>-1
\eeqnn
and $\gamma_{n,ij}(\d s)\geq \bar{b}_{ij}(\d s)$. Also we have
\beqnn
	\int_0^t\int_{\mbb{R}_+^2\setminus\{\mbf{0}\}}(z_1+z_2)\mu_{n,i}(\d s,\d\mbf{z})\ar=\ar2nc_i(t)+(1-\e^{-n})\int_0^t\int_{\mbb{R}_+^2\setminus\{\mbf{0}\}}(z_i+z_j)(1\wedge n\Vert \mbf{z}\Vert)m_i(\d s,\d\mbf{z})\\
	\ar\leq\ar 2nc_i(t)+(n+1)m_i(t)<\infty.
\eeqnn

\begin{lemma}\label{phi_n and phi properties}
	The branching mechanism $\bm\phi$ and $\bm\phi_n$ have the following properties:
	\begin{itemize}
		\item[\rm(1)] For $t\geq r\geq 0$ and $\mbf f\in\mbf B[0,\infty)^+$, we have $\phi_i(\mbf f,(r,t])=\uparrow\lim_{n\to\infty}\phi_{n,i}(\mbf f,(r,t])$;
		\item[\rm(2)] For $t\geq s\geq r\geq 0$ and $\mbf{f,g}\in\mbf B[0,\infty)^+$ satisfying $f_k\leq g_k,\ k=1,2$, we have
		\beqnn
			\phi_i(\mbf f,(s,t])-\phi_{n,i}(\mbf f,(s,t])\leq \phi_i(\mbf g,(r,t])-\phi_{n,i}(\mbf g,(r,t]);
		\eeqnn
		\item[\rm(3)] For $t\geq r\geq 0$ and $\mbf{f,g}\in\mbf B[0,\infty)^+$, we have
		\beqnn
		\begin{split}
			\sup_i|\phi_i(\mbf f,(r,t])-\phi_i(\mbf g,(r,t])|\leq C_1(t)\int_r^t\sup_i |f_i(u)-g_i(u)| C_2(\d u),
		\end{split}
		\eeqnn
		where
		\beqnn	
			C_1(t)=\sup_{s\in[0,t]}[f_1(s)+f_2(s)+g_1(s)+g_2(s)]+1,
		\eeqnn
		and
		\beqlb\label{C_2}
		\begin{aligned}
			C_2(\d u)&=c_1(\d u)+2\int_{\mbb{R}_+^2\setminus\{\mbf{0}\}}(z_1\mbf 1_{\{\Vert\mbf{z}\Vert>1\}}+z_1^2\mbf 1_{\{\Vert\mbf{z}\Vert\leq 1\}}+z_2)m_1(\d u,\d\mbf{z})\\
			&\quad+c_2(\d u)+2\int_{\mbb{R}_+^2\setminus\{\mbf{0}\}}(z_2\mbf 1_{\{\Vert\mbf{z}\Vert>1\}}+z_2^2\mbf 1_{\{\Vert\mbf{z}\Vert\leq 1\}}+z_1)m_2(\d u,\d\mbf{z})\\
			&\quad+\Vert b_{11}\Vert(\d u)+\Vert b_{22}\Vert(\d u)+b_{12}(\d u)+b_{21}(\d u).
		\end{aligned}			
		\eeqlb
	\end{itemize}
\end{lemma}
\bproof(1) and (2) are obvious by the definition. For $t\geq r\geq 0$ and $\mbf{f,g}\in\mbf B[0,\infty)^+$, we have
\beqnn
	\ar\ar|\phi_i(\mbf f,(r,t])-\phi_i(\mbf g,(r,t])|\leq\int_r^t |f_i(u)-g_i(u)|\Vert b_{ii}\Vert (\d u)+\int_r^t |f_j(u)-g_j(u)|b_{ij}(\d u)\\
	\ar\ar\qquad+C_1(t)\int_r^t|f_i(u)-g_i(u)|c_i(\d u)+C_1(t)\int_r^t\int_{\mbb{R}_+^2\setminus\{\mbf{0}\}}|f_j(u)-g_j(u)|z_jm_i(\d u,\d\mbf{z})\\
	\ar\ar\qquad+C_1(t)\int_r^t\int_{\mbb{R}_+^2\setminus\{\mbf{0}\}}|f_i(u)-g_i(u)|\big(z_i\mbf 1_{\{\Vert\mbf{z}\Vert>1\}}+z_i^2\mbf 1_{\{\Vert\mbf{z}\Vert\leq 1\}}+z_j\big)m_i(\d u,\d\mbf{z}),
\eeqnn
implying (3).
\eproof

\bproposition\label{uniqueness}
Let $(b_{ij}, b_{ii},c_i,m_i)$ be defined as in introduction. Then for $\bm\lambda\in \mbb{R}_+^2$ and $t\geq 0$, there is at most one solution to \eqref{backward eq}. Moreover, for $t\geq r\geq 0,\ \bm\lambda \in \mbb{R}_+^2$, we have
	\beqlb\label{upper bound}
		v_{i,r,t}(\bm\lambda)\leq U_{i}(r,t,\bm\lambda),
	\eeqlb
	where
	\beqnn
		U_{i}(r,t,\bm\lambda)=\Vert\bm\lambda\Vert(1+\bar{b}_{ij}(t))\exp\left\{\e^{\Vert b_{jj}\Vert(t)}\bar{b}_{12}(t)\bar{b}_{21}(t)+\Vert b_{11}\Vert(t)+\Vert b_{22}\Vert(t)\right\}
	\eeqnn
	and $\bar{b}_{ij}$ is defined by \eqref{b_ij-}.
\eproposition

\bproof Suppose that $r\mapsto \mbf{v}_{r,t}(\bm\lambda)$ is a bounded positive solution to \eqref{backward eq}, we can easily obtain
 \beqnn
	\begin{split}
		v_{i,r,t}(\bm\lambda)\leq \lambda_i+\int_r^t v_{i,s,t}(\bm\lambda) \Vert b_{ii}\Vert(\d s)+\int_r^t v_{j,s,t}(\bm\lambda)\bar{b}_{ij}(\d s).
	\end{split}
 \eeqnn
Then by Lemma \ref{multiGronwall}, we have the estimation \eqref{upper bound}. Suppose that $r\mapsto \mbf w_{r,t}(\bm\lambda)$ is also a bounded positive solution to \eqref{backward eq}. By \eqref{upper bound}, we have
	\beqnn
		\ar\ar|v_{i,r,t}(\bm\lambda)-w_{i,r,t}(\bm\lambda)|\leq \int_r^t|v_{i,s,t}(\bm\lambda)-w_{i,s,t}(\bm\lambda)|\Vert b_{ii}\Vert(\d s)\\
		\ar\ar\quad+\int_r^t|v_{j,s,t}(\bm\lambda)-w_{j,s,t}(\bm\lambda)| b_{ij}(\d s)+2U_i(0,t,\bm\lambda)\int_r^t|v_{i,s,t}(\bm\lambda)-w_{i,s,t}(\bm\lambda)|c_i(\d s)\\
		\ar\ar\quad+\int_r^t\int_{\mbb{R}_+^2\setminus{\bar{B}^+(\mbf{0},1)}}\big[2|v_{i,s,t}(\bm\lambda)-w_{i,s,t}(\bm\lambda)|z_i+|v_{j,s,t}(\bm\lambda)-w_{j,s,t}(\bm\lambda)|z_j\big]m_i(\d s,\d\mbf{z})\\
		\ar\ar\quad+\int_r^t\int_{\bar{B}^+(\mbf{0},1)\setminus\{\mbf{0}\}} \Big[|v_{i,s,t}(\bm\lambda)-w_{i,s,t}(\bm\lambda)|\big(U_i(0,t,\bm\lambda)z_i^2+U_j(0,t,\bm\lambda)z_j\big)\\
		\ar\ar\quad\quad\quad\quad\quad\quad\quad\quad+|v_{j,s,t}(\bm\lambda)-w_{j,s,t}(\bm\lambda)|\big(U_1(0,t,\bm\lambda)+U_2(0,t,\bm\lambda)+1\big)z_j\Big]m_i(\d s,\d\mbf{z}).
	\eeqnn
	By Lemma \ref{multiGronwall}, we have $|v_{i,r,t}(\bm\lambda)-w_{i,r,t}(\bm\lambda)|=0$ for $r\in[0,t]$, implying the uniqueness.\eproof

\begin{remark} If we consider a weaker condition
 \beqnn
\bar m_i(t)=\int_0^t\int_{\mbb{R}_+^2\setminus\{\mbf{0}\}}(\Vert\mbf{z}\Vert^2\wedge 1+z_j\mbf{1}_{\{\Vert\mbf{z}\Vert\leq 1\}})m_i(\d s,\d\mbf{z})<\infty,\quad t\geq 0,
 \eeqnn
the estimation in above proposition will be more complex since we need a uniformly strictly positive lower bound of solution on the exponent to control the linear item of integrand on $\mbb{R}_+^2\setminus{\bar{B}^+(\mbf{0},1)}$. \end{remark}

\noindent\textit{Proof of Theorem~\ref{main result}.~}
Let $\bm\phi$ be a branching mechanism with parameters $(\bar{b}_{ij},b_{ii},c_i,m_i)$, where $\bar{b}_{ij}$ is defined by \eqref{b_ij-}. It is obvious that
\beqlb\label{3 backward eq}
	\mbf{v}_{r,t}(\bm\lambda)=\bm\lambda-\bm\phi(\mbf{v}_{\cdot,t}(\bm\lambda),(r,t]),\quad r\in [0,t],
\eeqlb
is equivalent to \eqref{backward eq}. Let $\bm\phi_n$ be defined by \eqref{phi_n}. By above arguments, we see that $\gamma_{n,ii}$, $\gamma_{n,ij}$ and $\mu_{n,i}$ satisfies the conditions of Proposition \ref{special case backward eq solution}. Then we can define a cumulant semigroup $((\mbf{v}^{(n)}_{r,t}(\bm\lambda))_{t\geq r}$ solving
\beqlb
	\mbf{v}_{r,t}(\bm\lambda)=\bm\lambda-\bm\phi_n(\mbf{v}_{\cdot,t}(\bm\lambda),(r,t]),\quad r\in [0,t].
\eeqlb
Also we have the estimation $v_{i,r,t}^{(n)}(\bm\lambda)\leq 2A \e^{\rho(t)}$ for $r\in[0,t]$ and $\bm\lambda\in[0,A]^2$, where
\beqnn
	\rho(t)=2\Vert b_{11}\Vert(t)+2\Vert b_{22}\Vert(t)+b_{12}(t)+b_{21}(t).
\eeqnn
For $n\geq k\geq 1$, define
\beqnn
	D_{k,n}(r,t,\bm\lambda)=\sup_{i}\sup_{r\leq s\leq t}|v_{i,s,t}^{(n)}(\bm\lambda)-v_{i,s,t}^{(k)}(\bm\lambda)|.
\eeqnn
Then by Lemma \ref{phi_n and phi properties}, we can easily obtain
\beqnn
\begin{split}
	D_{k,n}(r,t,\bm\lambda)&\leq\tilde A(t)+C_1(t)\int_r^t D_{k,n}(s,t,\bm\lambda)C_2(\d u),
\end{split}
\eeqnn
where
\beqnn
	\tilde A(t)=2\sup_i|\phi_i(2A \e^{\rho(t)},2A \e^{\rho(t)},(0,t])-\phi_{k,i}(2A \e^{\rho(t)},2A \e^{\rho(t)},(0,t])|,
\eeqnn
$C_1(t)=8A \e^{\rho(t)}+1$ and $C_{2}(\d u)$ is defined by \eqref{C_2}. Moreover, by Gronwall's inequality, we have
\beqnn
	D_{k,n}(r,t,\bm\lambda)\leq \tilde A(t)\e^{C_1(t)C_{2}(t)}.
\eeqnn
Then we can easily obtain the limit $v_{i,r,t}(\bm\lambda)=\lim_{k\to\infty}v_{i,r,t}^{(k)}(\bm\lambda)$ exists and convergence is uniform in $(r,\bm\lambda)\in[0,t]\times[0,A]^2$ for every $A\geq 0$. By Dominated convergence and Lemma \ref{phi_n and phi properties}, we see that the limit $r\mapsto \mbf{v}_{r,t}(\bm\lambda)$ is the solution to \eqref{3 backward eq}. Moreover, by Lemma 1 in Watanabe (1969), we see that $v_{i,r,t}(\bm\lambda)$ has representation \eqref{LK representation}. Recall that uniqueness is obtained in Proposition \ref{uniqueness}, implying the semigroup property \eqref{semigroup property}. Thus $((\mbf{v}_{r,t}(\bm\lambda))_{t\geq r}$ is a cumulant semigroup. \qed

\noindent\textit{Proof of Theorem~\ref{main result(2)}.~} Step 1. The uniqueness follows by Lemma \ref{multiGronwall}. Next we only need to prove the case for $\bm\lambda\in\mbb R_+^2$. In fact, if \eqref{first moments} holds for $\bm\lambda\in\mbb R_+^2$ with $[0,t]\ni r\mapsto \bm{\pi}_{r,t}(\bm\lambda)\in \mbb{R}_+^2$ solving \eqref{linear backward eq}, then for $\bm\lambda=(\sgn\lambda_1|\lambda_1|,\sgn\lambda_2|\lambda_2|)\in\mbb R^2$,
 \beqnn
\int_{\mbb R_+^2}\langle\bm\lambda,\mbf y\rangle Q_{r,t}(\mbf x,\d \mbf y)
 \ar=\ar
\sgn\lambda_1\int_{\mbb R_+^2}\langle(|\lambda_1|,0),\mbf y\rangle Q_{r,t}(\mbf x,\d \mbf y)+\sgn\lambda_2\int_{\mbb R_+^2}\langle(0,|\lambda_2|),\mbf y\rangle Q_{r,t}(\mbf x,\d \mbf y)\\
 \ar=\ar
\sgn\lambda_1\langle\mbf x,\bm\pi_{r,t}((|\lambda_1|,0))\rangle+\sgn\lambda_2\langle\mbf x,\bm\pi_{r,t}((0,|\lambda_2|))\rangle,
\eeqnn
where $[0,t]\ni r\mapsto \bm{\pi}_{r,t}((|\lambda_1|,0))\in \mbb{R}_+^2$ and $[0,t]\ni r\mapsto \bm{\pi}_{r,t}((0,|\lambda_2|))\in \mbb{R}_+^2$ are solutions to \eqref{linear backward eq} with $\bm\lambda$ replacing by $(|\lambda_1|,0)$ and $(0,|\lambda_2|)$, respectively. Let $\bm\alpha_{r,t}(\bm\lambda)=\sgn\lambda_1\bm\pi_{r,t}((|\lambda_1|,0))+\sgn\lambda_2\bm\pi_{r,t}((0,|\lambda_2|))$. Then from \eqref{linear backward eq}, it follows that
\beqnn
	\alpha_{i,r,t}(\bm\lambda)\ar=\ar\lambda_i+\sgn\lambda_1\int_r^t\pi_{j,s,t}((|\lambda_1|,0))\bar b_{ij}(\d s)+\sgn\lambda_2\int_r^t\pi_{j,s,t}((0,|\lambda_2|))\bar b_{ij}(\d s)\\
	\ar\ar-\sgn\lambda_1\int_r^t \pi_{i,s,t}((|\lambda_1|,0)) b_{ii}(\d s)-\sgn\lambda_2\int_r^t\pi_{i,s,t}((0,|\lambda_2|)) b_{ii}(\d s)\\
	\ar=\ar\lambda_i+\int_r^t\alpha_{j,s,t}(\bm\lambda)\bar b_{ij}(\d s)-\int_r^t \alpha_{i,s,t}(\bm\lambda) b_{ii}(\d s),
\eeqnn
where we use $\sgn\lambda_1(|\lambda_1|,0)+\sgn\lambda_2(0,|\lambda_2|)=\bm\lambda$ for the first equality.

Step 2. For $\bm\lambda\in\mbb R_+^2, a>0$, by Proposition \ref{uniqueness}, we have $a^{-1}v_{i,r,t}(a\bm\lambda)\leq a^{-1}U_i(r,t,a\bm\lambda)=U_i(r,t,\bm\lambda)$. It is obvious that $\bm\lambda\mapsto \mbf v_{r,t}(\bm\lambda)$ is continuous, so we can define $\pi_{i,r,t}(\bm\lambda)=\frac{\partial}{\partial a}v_{i,r,t}(a\bm\lambda)\big|_{a=0}$. Then by differentiating both sides of \eqref{3 backward eq} with $\bm\lambda$ replacing by $a\bm\lambda$ and using bounded convergence theorem, we have
\beqnn
	\pi_{i,r,t}(\bm\lambda)=\lambda_i+\int_r^t \pi_{j,s,t}(\bm\lambda)\bar b_{ij}(\d s)-\int_r^t \pi_{i,s,t}(\bm\lambda)b_{ii}(\d s),
\eeqnn
where we use the fact that $\mbf v_{r,t}(\mbf 0)=\mbf 0$. Moreover, we can differentiate both sides of \eqref{Laplace transform} with $\bm\lambda$ replacing by $a\bm\lambda$ to obtain
\beqnn
	\int_{\mbb R_+^2}\langle\bm\lambda,\mbf y\rangle Q_{r,t}(\mbf x,\d\mbf y)=-\frac{\partial}{\partial a}\int_{\mbb{R}_+^2}\e^{-a\langle {\bm{\lambda}},\mbf{y}\rangle}Q_{r,t}(\mbf{x},\d\mbf{y})\Big|_{a=0}
 =
-\frac{\partial}{\partial a}\e^{-\langle \mbf{x},\mbf{v}_{r,t}(a{\bm\lambda})\rangle}\Big|_{a=0}=\langle \mbf{x},\bm\pi_{r,t}({\bm\lambda})\rangle.
\eeqnn
\qed

\vskip 0.3cm {\bf Acknowledgements}
We are grateful to the Laboratory of Mathematics and Complex Systems (Ministry of Education) for providing us the research facilities. This research is supported by the National Key R\&D Program of China (No. 2020YFA0712901).

\vskip 0.3cm {\bf Conflict of Interest}\quad The authors declare no
conflict of interest.


\begin{thebibliography}{99}\small


\bibitem{AlS82} Aliev, S.A. and Shchurenkov, V.M. (1982). Transitional phenomena and the convergence of Galton--Watson processes to Ji\v{r}ina processes. \textit{Theory Probab. Appl.} \textbf{27}, 472--485.


\bibitem{BaS15} Bansaye, V. and Simatos, F. (2015). On the scaling limits of Galton--Watson processes in varying environment. \textit{Electron. J. Probab.} \textbf{20}, no.~75, 1--36.

\bibitem{FaL22} Fang, R. and Li, Z. (2022). Construction of continuous-state branching processes in varying environments. \textit{Ann. Appl. Probab.} \textbf{32}, 3645--3673.

\bibitem{Gri74} Grimvall, A. (1974). On the convergence of sequences of branching processes. \textit{Ann. Probab.} \textbf{2}, 1027--1045.

\bibitem{GHP24} Gumenyuk, P.; Hasebe, T. and P\'{e}rez, J.-L. (2024): Loewner Theory for Bernstein functions I: evolution families and differential equations. \textit{Constr. Approx.} Online: doi.org/10.1007/s00365-023-09675-9; \textit{arXiv:2206.04753}.

\bibitem{GHP22+} Gumenyuk, P.; Hasebe, T. and P\'{e}rez, J.-L. (2022+): Loewner Theory for Bernstein functions II: applications to inhomogeneous continuous-state branching processes. \textit{arXiv:2211.12442}.

\bibitem{Lam67} Lamperti, J. (1967). The limit of a sequence of branching processes. \textit{Z. Wahrsch. verw. Ge.} \textbf{7}, 271--288.

\bibitem{LiLi24} Li, P. and Li, Z. (2024). Uniqueness problem for the backward differential equation of a continuous-state branching process. \textit{Acta Math. Sin. (Engl. Ser.)} \textbf{40}, 1825--1836.

\bibitem{Li20} Li, Z. (2020). Continuous-state branching processes with immigration. In \textit{From Probability to Finance-Lecture Notes of BICMR Summer School on Financial Mathematics. Math. Lect. Peking Univ.} 1-69.

\bibitem{Mao90} Mao, X. (1990). Lebesgue-stieltjes integral inequalities in several variables with retardation. \textit{Proc. Indian Acad. Sci. (Math. Sci.)} \textbf{100}, 231-243.

\bibitem{RhS70} Rhyzhov, Y.M. and Skorokhod, A.V. (1970). Homogeneous branching processes with a finite number of types and continuous varying mass. \textit{Theory Probab. Appl.} \textbf{15}, 704--707.

\bibitem{Wat69} Watanabe, S. (1969). On two dimensional Markov processes with branching property. \textit{Trans. Amer. Math. Soc.} \textbf{136}, 447--466.


\end{thebibliography}
\end{document}